# On Ramsey Numbers

Dhananjay P. Mehendale
Sir Parashurambhau College, Tilak Road, Pune 411030,
India

## Abstract

In this paper we define new numbers called the Neo-Ramsay numbers. We show that these numbers are in fact equal to the Ramsay numbers. Neo-Ramsey numbers are easy to compute and for finding them it is not necessary to check all possible graphs but enough to check only special kind of graphs having a well-defined adjacency pattern.

1. **Introduction:** For the most comprehensive and most recent knowledge on small Ramsey numbers one should refer to [1]. The Ramsey number $r(k,l)$ stand for the smallest positive integer $r$ such that every graph on $r$ points contains either a complete graph on $k$ points or an independent set on $l$ points. Clearly, $r(k,l)-1$ is the **largest** integer such that we can construct a graph on $(r-1)$ points which neither contains a complete graph on $k$ points nor an independent set on $l$ points.
       We define the new number called the Neo-Ramsey number, $R(k,l)$, and show that $r(k,l) = R(k,l)$. We will offer a procedure to determine this new number $R(k,l)$. The generalized Ramsey number $r(m_1, m_2, \cdots, m_n)$ stand for the smallest positive integer $r$ such that every graph on $r$ points, whose edges are colored with some distinct $n$ colors, $c_1, c_2, \cdots, c_n$ say, contains at least one monochromatic complete graph on $m_i$ points with edges colored by color $c_i$, for some $i = 1, 2, \ldots, n$. The definition of Neo-Ramsey number can be easily generalized to Generalized-Neo-Ramsey number, $R(m_1, m_2, \cdots, m_n)$, and it can further be shown that this Generalized-Neo-Ramsey number, $R(m_1, m_2, \cdots, m_n)$, is related to standard generalized Ramsey number, $r(m_1, m_2, \cdots, m_n)$, by the following simple relation: $r(m_1, m_2, \cdots, m_n) = R(m_1, m_2, \cdots, m_n)$.
       Suppose M is the set of $p$ integers: M = {1, 2, ..., $p$}. Let $M_D$ be the difference set of M, i.e. the set arrived at by taking

mutual differences of distinct integers in M. Thus, $M_D$ = {1, 2, 3, ..., $p-1$}. Suppose there exists **at least one** disjoint partitioning of $M_D$ into two subsets A and B, i.e. $M_D$ = A $\cup$ B and A $\cap$ B = $\phi$, a null set, such that for every $(k-1)$-subset U of A and for every $(l-1)$-subset V of B the set $U^* \not\subset$ A and the set $V^* \not\subset$ B, where $U^*$ = U $\cup$ { set of integers not in U and arrived at by taking mutual differences among each other of all distinct elements in U}, and $V^*$ = V $\cup$ {set of integers not in V and arrived at by taking mutual differences among each other of all distinct elements in V}. This property satisfied by set A, B is called $(k,l)$- **property**.

The set U* will be called hereafter the **closure** of the set U (The set V* will be called hereafter the closure of the set V) and the operation of constructing U* from U (V* from V) will be described hereafter as **taking the closure** of U (V). If U* = U (V* = V) then the set U (V) will be called a **closed set**. Note that a set of numbers forming a progression is a closed set: e.g. if U = {7, 14, 21, 28} then U* = U.

We now proceed to define the Neo-Ramsey number, $R(k,l)$, in terms of the largest set M whose difference set $M_D$ has some disjoint subsets A, B whose union is $M_D$ and possessing $(k,l)$ property.

**Definition 1.1:** Let M be the **largest set** (with cardinality $p$) whose difference set $M_D$ with cardinality $p-1$ possesses the $(k,l)$ property then the **Neo-Ramsay number**, $R(k,l) = p+1$.

Note that by largest we meant that the difference set of any set N = {1, 2, ..., $q$}, such that $q > p$, does not possess $(k,l)$-property.

In order to understand the **motivation** behind defining the Neo-Ramsay number we proceed with the following alternative equivalent **tabular representation (bitableau)** of the adjacency matrix of the graph.

Let $G$ be a (p, q) graph, i.e. a graph on p points (vertices) and q lines (edges) with the following vertex set $V(G)$ and edge set $E(G)$ respectively:

$$V(G) = \{v_1, v_2, \cdots, v_p\} \text{ and}$$
$$E(G) = \{e_1, e_2, \cdots, e_q\}$$

**Definition 1.2:** The **increasing vertex adjacency bitableau** associated with every labeled copy of graph $G$, $IVAB(G)$, is the following bitableau:

$$IVAB(G) = \begin{pmatrix} 1 & 2 & 3* & \cdots & p \\ 2 & 3 & 4* & \cdots & p* \\ \vdots & & & & \\ j & (j+1)* & (j+2) & \cdots & p \\ \vdots & & & & \\ p & - & & \cdots & \end{pmatrix}$$

where left tableau represents the suffixes of the vertex labels and stand for the vertices while the right tableau represents the rows of the suffixes of the vertex labels and represent the vertices that are adjacent (nonadjacent) to the vertex whose suffix is written in the same row in the left tableau, i.e. the appearance of **entry $k\,(k*)$ in the $j^{th}$ row** of the right tableau implies that $j < k$ and vertex $v_j$ is **adjacent (nonadjacent)** to vertex $v_k$.

Thus, for example, from the above given $IVAB(G)$ vertex $v_1$ is adjacent to vertex $v_p$, but vertex $v_j$ is not adjacent to vertex $v_{(j+1)}$, etc. It is clear to see that for a complete graph, $K_p$, on $p$ points

$$IVAB(K_p) = \begin{pmatrix} 1 & 2 & 3 & \cdots & p \\ 2 & 3 & 4 & \cdots & p \\ \vdots & & & & \\ j & (j+1) & (j+2) & \cdots & p \\ \vdots & & & & \\ p & - & & \cdots & \end{pmatrix}$$

Similarly, for an independent set, $I_p$, on $p$ points

$$IVAB(I_p) = \begin{pmatrix} 1 & 2* & 3* & \cdots & p* \\ 2 & 3* & 4* & \cdots & p* \\ \vdots & & & & \\ j & (j+1)* & (j+2)* & \cdots & p* \\ \vdots & & & & \\ p & - & & \cdots & \end{pmatrix}$$

**Example 1.1:** Let $k = l = 3$. For this case, we know that $r(k,l) = r(3,3) = 6$. So, We take M = {1, 2, 3, 4, 5}. So, $M_D$ = {1, 2, 3, 4}. We now express $M_D$ as union of two 2-element disjoint subsets A and B.
(1) We take A = {1, 2} and B = {3, 4}. Since cardinality of A = cardinality of B = 2, therefore, there is only one $(k-1)$-subset = 2-subset U of A, namely, the set A itself, and there is only one $(l-1)$-subset = 2-subset V of B, namely, the set B itself. It is clear from the definition that U* = {1, 2} and V* = {1, 3, 4}, and therefore, U* ⊂ A and V* ⊄ B.
(2) We take A = {1, 3} and B = {2, 4}. Again, since cardinality of A = cardinality of B = 2, therefore, there is only one $(k-1)$-subset = 2-subset U of A, namely, the set A itself, and there is only one $(l-1)$-subset = 2-subset V of B, namely, the set B itself. It is clear from the definition that U* = {1, 2, 3} and V* = {2, 4}, and therefore, U* ⊄ A and V* ⊂ B.
(3) We take A = {1, 4} and B = {2, 3} such that A∩B = $\phi$. Since cardinality of A = cardinality of B = 2, therefore, there is only one $(k-1)$-subset = 2-subset U of A, namely, the set A itself, and there is only one $(l-1)$-subset = 2-subset V of B, namely, the set B itself. It is clear from the definition that U* = {1, 3, 4} and V* = {1, 2, 3}, and therefore, U* ⊄ A and V* ⊄ B.

Now, let us define a graph G, on five points {$v_1, v_2, \cdots, v_5$}, such that points ($v_i, v_j$) are adjacent if the difference $|i - j| \; \varepsilon$ A and nonadjacent if the difference $|i - j| \; \varepsilon$ B. It is easy to check for the case (1) and the case (2) that the graph will contain either a complete graph on three points or an independent set on three points. But, for the case (3) the graph will neither contain a complete graph on three points nor an independent set on three points! The $IVAB(G)$ for this third case, case (3), will be

$$IVAB(G) = \begin{pmatrix} 1 & 2 & 3* & 4* & 5 \\ 2 & 3 & 4* & 5* & \\ 3 & 4 & 5* & & \\ 4 & 5 & & & \\ 5 & - & & & \end{pmatrix}$$

Now, if one works out this procedure with M = {1, 2, 3, 4, 5, 6}, so that, $M_D$ = {1, 2, 3, 4, 5}, and define all possible A sets containing 3 (or 2) elements of $M_D$ and B sets containing (remaining) 2 (or 3) elements of $M_D$, or vice versa, and consider all possible 2-subsets U and V of each A and corresponding B we can see that **not a single** disjoint partitioning of $M_D$ will satisfy $U^* \not\subset A$ and $V^* \not\subset B$. Thus, 5 is the **largest** number satisfying the (3, 3)-property so, $R(3,3) = 6$, which is actually equal to $r(3,3)$!

**2: The relationship between $r(k,l)$ and $R(k,l)$:** Note that for us the $(k-1)$-subset $U = \{j_1, j_2, \cdots, j_{(k-1)}\}$ of set A (or $(l-1)$-subset $V = \{j_1, j_2, \cdots, j_{(l-1)}\}$ of set B) represent the graphs with vertex sets $\{v_i, v_{(i+j_1)}, v_{(i+j_2)}, \cdots, v_{(i+j_{(k-1)})}\}$ (or $\{v_i, v_{(i+j_1)}, v_{(i+j_2)}, \cdots, v_{(i+j_{(l-1)})}\}$).

If $U^* \subset A$ (or $V^* \subset B$) then all these graphs are complete graphs (or independent sets) on $k$ (on $l$) points.

We now proceed to show that

**Lemma 2.1:** $r(k,l) > R(k,l) - 1$.

**Proof:** Let $M = \{1, 2, \ldots, p\}$ be the largest set possessing $(k,l)$ property. We proceed to show that we can construct a graph on $p$ points which neither contains a complete graph on $k$ points nor an independent set on $l$ points. We construct the graph G with desired property as follows:

Let $V(G) = \{v_1, v_2, \cdots, v_p\}$, be the vertex set of $G$. $M_D$ = {1, 2, 3, …, $p-1$}, the set of the mutual differences of the elements of M. Let $M_D = A \cup B$ and $A \cap B = \phi$ such that for every $(k-1)$-subset U of A and for every $(l-1)$-subset V of B the set $U^* \not\subset A$ and the set $V^* \not\subset B$, where $U^* = U \cup \{$ set of integers not in U and arrived at by taking mutual differences among each other of all distinct elements in U$\}$, and $V^* = V \cup \{$set of integers not in V and arrived at by taking mutual differences among each other of all distinct elements in V$\}$.

We define the adjacency of $G$ as follows:

(i) The vertices $v_i$, $v_j$ of $G$ are adjacent if the difference $|i - j| \; \varepsilon \; A$

(ii) The vertices $v_i$, $v_j$ of $G$ are nonadjacent if the difference $|i - j| \; \varepsilon \; B$

It is clear to see that the graph $G$ will neither contain a complete graph on $k$ points nor an independent set on $l$ points, as follows:

a) Consider a subgraph with vertex set $\{v_{i_1}, v_{i_2}, \cdots, v_{i_k}\}$. As mentioned above, associate with this graph the difference set, the $(k-1)$-subset $U = \{i_2 - i_1, i_3 - i_1, \cdots, i_k - i_1\}$. $U^* = U \cup \{$ numbers of type $(i_r - 1) - (i_s - 1) = i_r - i_s$ such that $r > s \}$. Since, $U^* \not\subset A$, therefore, the vertex set $\{v_{i_1}, v_{i_2}, \cdots, v_{i_k}\}$ will not form a complete graph on $k$ points.

b) Consider a subgraph with vertex set $\{v_{i_1}, v_{i_2}, \cdots, v_{i_l}\}$. As mentioned above, associate with this graph the difference set, the $(l-1)$-subset $V = \{i_2 - i_1, i_3 - i_1, \cdots, i_l - i_1\}$. $V^* = V \cup \{$ numbers of type $(i_r - 1) - (i_s - 1) = i_r - i_s$ such that $r > s \}$. Since, $V^* \not\subset B$, therefore, the vertex set $\{v_{i_1}, v_{i_2}, \cdots, v_{i_l}\}$ will not form an independent set on $l$ points. Hence etc. □

So far we have considered the disjoint partitioning of set $M_D$ as $M_D = A \cup B$ such that $A \cap B = \phi$. We hereafter call such a (disjoint) partitioning the **total or strict partitioning**. In the total partitioning we make vertices which are separated from each other by a particular fixed distance either adjacent (and so the number representing that distance is taken to be belonging to set A) or nonadjacent (and so the number representing that distance is taken to be belonging to set B). But instead if we take some vertex pairs which are separated from each other by a particular fixed distance adjacent while some other separated from each other by the same fixed distance nonadjacent then we call such a partitioning the **partial partitioning**.

To get a fill of the advantage of total partitioning we take a look at the following well-known theorem due to Turan [2]:

**Theorem 2.1 (Turan, 1941):**

$$ex(p, K_n) = \frac{p(p-1)}{2} - \{rp - \frac{(n-1)r(r+1)}{2}\}$$

where $r = \left[\frac{p-1}{n-1}\right]$, and $[x]$ = integral part of $x$. In this formula the first term represents the number of lines in the complete graph on $p$ points and the second term represents the (minimum) number of lines to be removed from this complete graph so that every complete graph on $n$ points gets eliminated.

**Proof:** We first construct a graph on $p$ points containing the number of lines equal to the right hand side of the above formula and not containing any clique on $n$ points. For establishing the maximality of this number one can carry out an easy check of the numerical equivalence of the above formula and the formula originally proposed by Turan [2].

Construct $IVAB(K_p)$ for complete graph on $p$ points. Mark the entries by a * (indicating the removal of corresponding lines) in the entire columns containing entry $n, 2n-1, 3n-2, 4n-3, \cdots$ thus, numbers starting with $n$ and are at a distance $n-1$ in succession, in the first row of the right tableau of $IVAB(K_p)$.

We can see that thus resulted increasing vertex adjacency bitableau does not contain any clique on $n$ points due to the following simple fact: Given any n numbers $\{j_1, j_2, \cdots, j_n\}$ then in the system of differences $|j_r - j_s|$, for all $r, s$ such that $r = 2, 3, \ldots, n$ and $s = 1, 2, \ldots, r-1,$ there exists at least one difference such that $|j_r - j_s| \equiv 0 \mod (n-1)$. Hence etc. □

**Illustration:** The $IVAB(G)$ constructed as described above for a graph on $p = 7$ points without containing a clique on $n = 4$ points:

$$IVAB(G) = \begin{pmatrix} 1 & 2 & 3 & 4* & 5 & 6 & 7* \\ 2 & 3 & 4 & 5* & 6 & 7 & \\ 3 & 4 & 5 & 6* & 7 & & \\ 4 & 5 & 6 & 7* & & & \\ 5 & 6 & 7 & & & & \\ 6 & 7 & & & & & \\ 7 & - & & & & & \end{pmatrix}$$

**Definition 2.1:** A partitioning of $M_D = A \cup B$ is called **partial partitioning** if some vertex pairs represented by a particular separation between the suffixes of the vertex labels are adjacent and so the number representing that distance belongs to A while some other vertex pairs represented by the same distance are nonadjacent and so the number representing that distance belongs to B also, and therefore, $A \cap B \neq \phi$.

A **better representation** for partial partitioning in terms of sets is as follows: Consider a graph $G$ on $p$ points, and let $V(G) = \{v_1, v_2, \cdots, v_p\}$, be the vertex set of $G$. Construct $IVAB(G)$.

(a) Read first **column** in the right tableau of $IVAB(G)$.
   (1) If it contains all entries **without** a * mark in front of it then put number 1 in set A.
   (2) If it contains all entries **with** a * mark in front of it then put number 1 in set B.
   (3) If it contains some entries **with** a * mark in front of them while some other entries **without** a * mark in front of them then put number 1 in set C.

(b) Read each column (from first to last column) and carry out the above procedure carried out for the first column for all columns and put numbers 2, 3, ..., $(p-1)$) in the appropriate sets, A, B, C.

Thus, $M_D = A \cup B \cup C$, $A \cap B = \phi$, and C represents the set of distances which represent adjacent vertex pairs for some vertices considered in pairs and nonadjacent vertex pairs for some other vertices considered in pairs for the same separation.

Thus, a partial partitioning can described by redefining sets A, B as $A = A \cup C$, $B = B \cup C$ and $A \cap B = C \neq \phi$.

Our **objective** is to find the largest number for which we can show the **existence of a graph** which neither contains a complete graph on $k$ points nor an independent set on $l$ points. This number is largest even if we allow total or partial partitioning. By definition we have denoted such a largest number for total partitioning as $R(k,l)-1$. But to take into account **all possible graphs** and to show that in fact $R(k,l)-1$ is that largest number for all possible graphs having total or partial partitioning which neither contain a complete graph on $k$ points nor an independent set on $l$ points, we need prove the following :

**Lemma 2.2:** Partial partitioning is not better than the total partitioning, i.e. we cannot improve upon the **largest size** of the vertex set for a graph, $p = (R(k,l)-1)$, by considering partial partitioning instead of total partitioning, such that the new graph with vertices more than $p$ contains neither a complete graph on $k$ points nor an independent set on $l$ points.

**Proof:** Let $M = \{1, 2, ..., R(k,l)-1\}$, so that $M_D = \{1, 2, ..., R(k,l)-2\}$
Let $M_D = A \cup B$, such that $A \cap B = \phi$, where $B = \{\alpha_1, \alpha_2, \cdots, \alpha_m\}$

Using sets A, B one can construct a graph on $R(k,l)-1$ points whose $IVAB(G)$ is as follows:

$$IVAB(G) = \begin{pmatrix} 1 & 2 & \cdots & (\alpha_1+1)* & \cdots & (\alpha_m+1)* & \cdots & p \\ 2 & 3 & \cdots & (\alpha_1+2)* & \cdots & (\alpha_m+2)* & \cdots & \\ \vdots & \vdots & & \vdots & & \vdots & & \\ p-1 & p & & & & & & \\ p & - & & & & & & \end{pmatrix}$$

and this graph is of largest size for total partitioning which is free from any clique on $k$ points or an independent set on $l$ points (We call such a total partitioning a **legal partitioning** and such a graph a **legal graph**), and by largest we meant that we cannot add number $R(k,l)-1$ in set A or B without violation, namely, without the creation of a clique on $k$ points or an independent set on $l$ points. Now, the question is whether we can do so by some different partial partitioning? We now proceed to see that we cannot.

Note that every partial partitioning leads to **increase** in the size of sets A, B due to addition of the elements which are (partially) common to both A and B. As mentioned above, any partial partitioning can be expressed as A = A ∪ C, B = B ∪ C and A ∩ B = C ≠ ϕ. Thus, if we cannot avoid the existence of cliques(s) or independent set(s) in total partitioning we cannot do so in partial partitioning.

So, only those graphs (having a partial partitioning) which will be isomorphic labeled copies of graphs like $G$ given above obtained through a legal total partitioning will be the suitable graphs. Hence etc. □

**Theorem 2.2:** $R(k,l) = r(k,l)$.

**Proof:** From lemmas 2.1 and 2.2 the result follows. □

**3. An Algorithm to Find Sets A, B:** Since the sets A and their complements in $M_D$, namely, sets B are too many in number it becomes a difficult task to check all possible sets A and B for checking the $(k,l)$ property. So, instead of this approach it is wiser to see which sets A and B could be ignored as **illegal** (producing illegal partitioning of the set $M_D$). It is wiser to eliminate useless possibilities for A and B (e.g. the sets A, B having subsets of numbers of size $(k-1)$ or of size $(l-1)$ respectively forming a progression are inappropriate ones since they are closed sets) before checking $(k,l)$ property. So, instead of checking all possible sets A and B we dynamically construct these sets checking $(k,l)$ property **at each intermediate stage**. For this to achieve, we proceed to discuss a simple and algorithm to construct the appropriate sets A, B with desired properties.

**Definition 3.1:** A number *m* is called the **limiting number** if the sets A and B together contain numbers up to *m*, i.e. $M_D = \{1, 2, ..., m\}$, and the **sets satisfy $(k,l)$ property.**

### 3.1 Steps of the Algorithm:
(1) Start with the **limiting number** $m = 0$, so that, $A = \phi$, $B = \phi$.
(2) Increase the limiting number by one, i.e. replace the limiting number *m* by $(m+1)$ and write down all A, B satisfying $(k,l)$ property, so that, $A = \{1\}$, $B = \phi$ (i.e. determine all possible disjoint partitionings satisfying $(k,l)$ property for the set containing numbers up to the present limiting number.)
(3) Go to (2), so that, (i) $A = \{1, 2\}$, $B = \phi$. (ii) $A = \{1\}$, $B = \{2\}$, if $(k,l)$ property is satisfied for both the possibilities (i) and (ii), and discard those possibilities for which $(k,l)$ property is not satisfied. (i.e. determine all possible disjoint partitionings satisfying $(k,l)$ property for the set containing numbers up to the present limiting number.)
(4) Continue adding next *m* to sets A, B in the all legal (satisfying $(k,l)$ property for the present limiting number) combinations of sets A, B till that *m* gets placed in some set A or B of some earlier legal combination producing at least one legal combination for the next stage, and stop at the stage $m \,(= R(k,l)-2)$ for which we cannot have any one legal combination declaring $R(k,l)$ as Ramsey number.

We now proceed, for the sake of **illustration**, to apply the algorithm for the simplest case of $r(3,3)$ as follows:
   (i)   Let $m = 0$, therefore, $A = \phi$, $B = \phi$.
   (ii)  Let (next *m*) = (earlier *m*) +1 = 1, therefore, $A = \{1\}$, $B = \phi$.
   (iii) Let (next *m*) = (earlier *m*) +1 = 2, so, $A = \{1\}$, $B = \{2\}$.
   (iv)  Let (next *m*) = (earlier *m*) +1 = 3, so, $A = \{1\}$, $B = \{2, 3\}$.
   (v)   Let (next *m*) = (earlier *m*) +1 = 4, so, $A = \{1, 4\}$, $B = \{2, 3\}$.
   (vi)  Let (next *m*) = (earlier *m*) +1 = 5, It is clear that this case we cannot extend the A, B sets at the earlier stage (step(v)) to produce extended new valid A, B sets corresponding to the present stage. Hence, $r(3,3) = 6$.

**More Examples:**
(1) We know that $r(3,4) = 9$. By using the above algorithm, we get at last (i.e. non-extendable) stage A = {1, 4, 7} and B = {2, 3, 5, 6}, or, A = {3, 4, 5} and B = {1, 2, 6, 7}.
(2) We know that $r(3,5) = 14$. By using the above algorithm, we get at last stage A = {1, 5, 8, 12} and B = {2, 3, 4, 6, 7, 9, 10, 11}, or A = {4, 6, 7, 9} and B = {1, 2, 3, 5, 8, 10, 11, 12}.
(3) We know that $r(3,9) = 36$. By using the above algorithm, we get A = {8, 12, 14, 17, 18, 21, 23, 27} and B = {1, 2, 3, 4, 5, 6, 7, 9, 10, 11, 13, 15, 16, 19, 20, 22, 24, 25, 26, 28, 29, 30, 31, 32, 33, 34}.
(4) We know that $r(4,4) = 18$. By using the above algorithm, we get A = {1, 2, 4, 8, 9, 13, 15, 16} and B = {3, 5, 6, 7, 10, 11, 12, 14}.
(5) We know that $r(4,5) = 25$. By using the above algorithm, we get A = {4, 7, 8, 9, 10, 14, 15, 16, 17, 20} and B = {1, 2, 3, 5, 6, 11, 12, 13, 18, 19, 21, 22, 23}.

**4. Generalized Neo-Ramsey numbers:** In order to determine the generalized Neo-Ramsey numbers, $R(m_1, m_2, \cdots, m_n)$, we will need to determine disjoint partitioning of the largest possible set $M_D$ into $n$ subsets, $A_1, A_2, A_3, \ldots, A_n$ such that, $M_D = \bigcup A_i$ and $\bigcap A_i = \phi$, $i = 1, 2, 3, \ldots, n$, and for every $(m_i - 1)$-subset $U_i$ of $A_i$ $U_i^* \not\subset A_i$ where $U_i^*$ is closure of $U_i$. And so, $R(m_1, m_2, \cdots, m_n) = \#(M_D) + 2$, since $M_D$ is the largest set. Note that everything discussed above for bicolored graphs can be straightforwardly extended for this general case.

**5. Conclusion:** The range of the known Ramsey numbers and even that of the multicolored Ramsey numbers can be made much higher with the help of the new approach discussed in this paper.

Acknowledgements

The author is very thankful to Dr. M. R. Modak, Bhaskaracharya Pratishthana, Pune, for many useful discussions, and Mr. Abhijit Patwardhan, Pune, for preparing a nice C program for the algorithm.